\documentclass[journal]{IEEEtran}

\usepackage{amsmath,amssymb,float,arydshln,color}
\usepackage{psfrag,setspace,wrapfig,subfigure}
\usepackage[latin1]{inputenc}
\usepackage{dsfont}
\usepackage[lined,ruled,commentsnumbered]{algorithm2e}
\usepackage{epsfig}
\usepackage{epstopdf}
\usepackage{graphicx}

\usepackage{amsfonts}
\usepackage{cite}
\usepackage{url}

\usepackage{hhline}
\usepackage[table]{xcolor}
\usepackage{multirow}
\usepackage{tabu}
\allowdisplaybreaks

\newcommand{\swb}{{\scriptstyle{\boldsymbol{\mathcal{W}}}}}
\newcommand{\ssb}{{\scriptstyle{\boldsymbol{\mathcal{S}}}}}

\newcommand{\sgb}{{{\boldsymbol{\mathcal{G}}}}}
\DeclareMathOperator*{\argmin}{arg\,min}
\newcommand{\bp}{ \begin{proof}}
\newcommand{\ep}{\end{proof} }

\newcommand{\Ex}{\mathbb{E}\hspace{0.05cm}}

\newcommand{\bm}[1]{\mbox{\boldmath $#1$}}

\newcommand{\be}{\begin{equation}}
\newcommand{\ee}{\end{equation}}
\newcommand{\bal}{\begin{align}}
	\newcommand{\eal}{\end{align}}
\newcommand{\bq}{\begin{eqnarray}}
\newcommand{\eq}{\end{eqnarray}}
\newcommand{\bqn}{\begin{eqnarray*}}
	\newcommand{\eqn}{\end{eqnarray*}}
\newcommand{\nn}{\nonumber}
\newcommand{\ba}{\left[ \begin{array}}
	\newcommand{\ea}{\\ \end{array} \right]}

\newcommand{\define}{\;\stackrel{\Delta}{=}\;}

\newtheorem{theorem}{{Theorem}}
\newtheorem{corollary}{{Corollary}}

\def\tran{^{\mathsf{T}}}

\def\one{\mathds{1}}

\def\bpsi	{{\boldsymbol \psi}}

\def\bgamma {{\boldsymbol \gamma}}

\def\h{{\boldsymbol{h}}}

\def\n{{\boldsymbol{n}}}

\def\s{{\boldsymbol{s}}}

\def\w{{\boldsymbol{w}}}
\def\x{{\boldsymbol{x}}}

\def\cJ{{\mathcal{J}}}

\def\cU{{\mathcal{U}}}
\def\cV{{\mathcal{V}}}

\def\real{{\mathbb{R}}}

\makeatletter
\def\hlinewd#1{%
	\noalign{\ifnum0=`}\fi\hrule \@height #1 \futurelet
	\reserved@a\@xhline}
\makeatother

\begin{document}
	\title{Performance Limits of Stochastic Sub-Gradient  Learning,\\ Part II: Multi-Agent Case}
	\author{Bicheng~Ying,~\IEEEmembership{Student Member,~IEEE},~and Ali~H.~Sayed,~\IEEEmembership{Fellow,~IEEE}
		\thanks{This work was supported in part by NSF grants CIF-1524250, ECCS-1407712, and DARPA N66001-14-2-4029. A short conference version appears in \cite{ying2016performance}.
			
			The authors are with Department of Electrical Engineering, University of California, Los Angeles, CA 90095. Emails: \{ybc,sayed\}@ucla.edu.}}
	
	\maketitle
	
\begin{abstract}
		The analysis in Part I\cite{ying16ssgd1} revealed interesting properties for subgradient learning algorithms in the context of {\em stochastic} optimization when gradient noise is present. These algorithms are used when the risk functions are non-smooth and involve non-differentiable components. They have been long recognized as being slow converging methods.
		However, it was revealed in Part I \cite{ying16ssgd1} that the rate of convergence becomes linear for {\em  stochastic} optimization problems, with the error iterate converging at an exponential rate $\alpha^i$ to within an $O(\mu)-$neighborhood of the optimizer, for some $\alpha \in (0,1)$ and small step-size $\mu$. The conclusion was established  under weaker assumptions than the prior literature and, moreover, several important problems (such as LASSO, SVM, and Total Variation) were shown to satisfy these weaker assumptions automatically (but not the previously used conditions from the literature). These results revealed that sub-gradient learning methods have more favorable behavior than originally thought when used to enable continuous adaptation and learning. The results of Part I\cite{ying16ssgd1} were exclusive to single-agent adaptation. The purpose of the current Part II is to examine the implications of these discoveries when a collection of networked agents employs subgradient learning as their cooperative mechanism. The analysis will show that, despite the coupled dynamics that arises in a networked scenario, the agents are still able to attain linear convergence in the stochastic case; they are also able to reach agreement within $O(\mu)$ of the optimizer.
	\end{abstract}
	

	\section{Introduction and Review of \cite{ying16ssgd1}}%
	We review briefly the notation and findings from Part I \cite{ying16ssgd1} in preparation for examining the challenges that arise in the multi-agent scenario. In Part I\cite{ying16ssgd1}, we considered an optimization problem of the form:
	\be w^{\star}\;=\;\arg\min_{w\in\real^{M}}\;J(w)\ee
		where the possibly {\em non-differentiable} but strongly-convex risk function $J(w)$ was expressed as the expectation of some convex but also possibly non-differentiable loss function $Q(\cdot)$, namely,
		\be
		J(w)\define \Ex\;Q(w;\x)
		\ee
		Here, the letter $\x$ represents the random data and the expectation
		operation is over the distribution of this data. The following sub-gradient algorithm was introduced and studied in Part I\cite{ying16ssgd1} for seeking $w^{\star}$:
		\bq
		\w_i&=& \w_{i-1} -\mu \widehat{g}(\w_{i-1})\label{eq.ssgd1}\\
		S_i&=&\kappa S_{i-1} + 1\label{eq.ssgd2}\\
		\bar{\w}_i&=&\left(1-\frac{1}{S_i}\right)\bar{\w}_{i-1}\;+\;\frac{1}{S_i} \w_i\label{eq.ssgd3}
		\eq
		with initial conditions $S_{0}=1$, $\w_{0}=0$, and $\bar{\w}_0=0$. Boldface notation is used for $\w_i$ to highlight its stochastic nature since the successive iterates are generated by relying on streaming data realizations for $\x$. Moreover, the scalar $\kappa\in[\alpha,1)$, where $\alpha=1-O(\mu)$ is a number close to one. The term $\widehat{g}(\w_{i-1})$ in \cite{polyak1987introduction} is an approximate sub-gradient at location $\w_{i-1}$; it is computed from the data available at time $i$ and approximates a true sub-gradient denoted by $g(\w_{i-1})$. This true sub-gradient is unavailable since $J(w)$ itself is unavailable in the stochastic context. This is because the distribution of the data $\x$ is unknown beforehand, which means that the expected loss function cannot be evaluated. The difference between a true sub-gradient vector and its approximation is gradient noise and is denoted by
		\be
		\s_i(\w_{i-1})\define \widehat{g}(\w_{i-1})-g(\w_{i-1})
		\ee
		\subsection{Data Model and Assumptions}
		\noindent The following three assumptions were motivated in Part I\cite{ying16ssgd1}:
		\begin{enumerate}
			\item[1.] $J(w)$ is $\eta-$strongly-convex so that $w^{\star}$ is unique. 	The strong convexity of $J(w)$ means that 
			\begin{align}
			J(\theta w_1+(1-\theta)w_2) \;\leq&\; \theta J(w_1) + (1-\theta) J(w_2)  \nn\\
			&\hspace{2mm}- \frac{\eta}{2} \theta(1-\theta)\|w_1-w_2\|^2,
			\end{align}
			for any $\theta\in[0,1]$, $w_1,$ and $w_2$. The above condition is equivalent to requiring \cite{polyak1987introduction}:
			\be
			J(w_1)\geq J(w_2) + g(w_2)\tran(w_1-w_2)+ \frac{\eta}{2} \|w_1-w_2\|^2 \label{assump.strongCVX}.
			\ee\label{sfewju}
			\item[2.] The subgradient is affine Lipschitz, meaning that there exist constants $c\geq 0$ and $d\geq 0$ such that 
			\be
			\|g(w_1)-g'(w_2)\|\;\leq\;c\|w_1-w_2\|\;+\;d,\;\;\;\forall w_1,w_2\label{assume.2}
			\ee
			and for any $g'(\cdot)\in\partial J(\cdot)$. Here, the notation $\partial J(w)$ denotes the differential at location $w$ (i.e., the set of all possible subgradient vectors at $w$). It was explained in Part I \cite{ying16ssgd1} how this affine Lipschitz condition is weaker than conditions used before in the literature and how important cases of interest (such as SVM, LASSO, Total Variation) satisfy it automatically (but do not satisfy the previous conditions). For later use, it is easy to verify (as was done in (50) in Part I\cite{ying16ssgd1}) that condition (\ref{assume.2}) implies that
			\be
				\|g(w_1) - {\color{black}g'(w_2)}\|^2 \leq e^2\|w_1 -w_2\|^2 + f^2, \;\; \forall w_1,w_2\label{assump.sbg2},
			\ee
			for any $g'(\cdot)\in\partial J(\cdot)$ and some constants $e^2\geq0$ and $f^2\geq0$. 
			\item[3.] The first and second-order moments of the gradient noise process satisfy the conditions:
			\begin{align}
					\Ex [\,s_{i}(\w_{i-1}) \,|\,\bm{\cal F}_{i-1}\,] =&\, 0\label{usadkh.13lk1l3k2},\\
					\hspace{-1mm}\Ex [\,\|s_{i}(\w_{i-1})\|^2 \,|\,\bm{\cal F}_{i-1}\, ] \leq&\, \beta^2 \|w^{\star}-\w_{i-1}\|^2 +\sigma^2\label{usadkh.13lk1l3k},
			\end{align}
			for some constants $\beta^2\geq0$ and $\sigma^2\geq 0$, and where the notation $\bm{\cal F}_{i-1}$ denotes the filtration (collection) corresponding to all past iterates:
			\be \bm{\cal F}_{i-1}\;=\;\mbox{\rm filtration by $\{\w_{j},\;j\leq i-1\}$}.
			\ee
			It was again shown in Part I\cite{ying16ssgd1} how the gradient noise process in important applications (e.g., SVM,LASSO) satisfy (\ref{usadkh.13lk1l3k2})---(\ref{usadkh.13lk1l3k}) directly. 
		\end{enumerate}
		Under the three conditions 1) --- 3), which are automatically satisfied for important cases of interest, the following important conclusion was proven in Part I \cite{ying16ssgd1} for the stochastic  subgradient  algorithm (\ref{eq.ssgd1})--(\ref{eq.ssgd3}) above. At every iteration $i$, it will hold that
		\be
		\lim_{i\to\infty}\Ex J(\bar\w_i) -J(w^\star) \leq  \mu(f^2+\sigma^2)/2 \label{eq.ff239ij3}
		\ee
		where the convergence of $\Ex J(\bar{\w}_i)$ to $J(w^{\star})$ occurs at an exponential rate $O(\alpha^i)$ where $\alpha=1-\mu\eta+O(\mu^2)$.
	
	\subsection{Interpretation of Result}
	For the benefit of the reader, we repeat here the interpretation that was given in Sec. IV.D of Part I \cite{ying16ssgd1} for the key results (\ref{eq.ff239ij3}); these remarks will be relevant in the networked case and are therefore useful to highlight again:
	
	\begin{enumerate}
		\item First, it has been observed in the optimization literature\cite{bertsekas1999nonlinear,polyak1987introduction,nesterov2004introductory} that sub-gradient descent iterations can perform poorly
		in deterministic problems (where $J(w)$ is known). Their convergence rate is $O(1/\sqrt{i})$
		under convexity and $O(1/i)$ under strong-convexity  when decaying step-sizes,
		$\mu(i)=1/i$, are used to ensure convergence \cite{nesterov2004introductory}. Result (\ref{eq.ff239ij3}) shows that the situation is different in
		the context of stochastic optimization when true subgradients
		are approximated from streaming data due to different requirements. By
		using {\em constant} step-sizes to enable continuous learning
		and adaptation, the sub-gradient iteration is now able to
		achieve exponential convergence at the rate of $O(\alpha^i)$ to steady-state.
		
		\item Second, this substantial improvement in convergence
		rate comes at a cost, but one that is acceptable
		and controllable. Specifically, we cannot guarantee convergence
		of the algorithm to the global minimum value,
		$J(w^{\star})$, anymore but can instead approach this optimal
		value with high accuracy in the order of $O(\mu)$, where
		the size of $\mu$ is under the designer's control and can be
		selected as small as desired.
		
		\item Third, this performance level is sufficient in most cases
		of interest because, in practice, one rarely has an infinite
		amount of data and, moreover, the data is often subject
		to distortions not captured by any assumed models. It
		is increasingly recognized in the literature that it is not
		always necessary to ensure exact convergence towards
		the optimal solution, $w^{\star}$, or the minimum value, $J(w^{\star})$,
		because these optimal values may not reflect accurately
		the true state due to modeling errors. For example, it is
		explained in the works \cite{bousquet2008tradeoffs,bottou2012stochastic,polyak1987introduction,towfic2014stability} that it is generally
		unnecessary to reduce the error measures below the
		statistical error level that is present in the data.
		
		\end{enumerate}
	
	\subsection{This Work}
	The purpose of this work is to examine how these properties reveal themselves in the networked case when a multitude of interconnected agents cooperate to minimize an aggregate cost function that is not generally smooth. In this case, it is necessary to examine closely the effect of the coupled dynamics and whether agents will still be able to agree fast enough under non-differentiability. 
	
	Distributed learning under non-smooth risk functions is common in many applications including distributed estimation and distributed machine learning. For example, $\ell_1$-regularization or hinge-loss functions (as in SVM implementations) lead to non-smooth risks.
	Several useful techniques have been developed in the literature for the solution of such distributed optimization problems, including the use of consensus strategies \cite{nedic2009distributed,yu2009distributed,kar2009distributed} and diffusion strategies \cite{sayed2014adaptive,chen2015learning1,chen2015learning2,tu2012diffusion}. In this paper, we will focus on the Adapt-then-Combine (ATC) diffusion strategy mainly because diffusion strategies have been shown to have superior mean-square-error and stability performance in adaptive scenarios where agents are expected to continually learn from streaming data\cite{tu2012diffusion}. 
	In particular, we shall examine the performance and stability behavior of networked diffusion learning under weaker conditions than previously considered in the literature. It is true that there have been several useful studies that employed sub-gradient constructions in the distributed setting before, most notably\cite{nemirovski2009robust,nedic2009distributed,ram2010distributed}. However, these earlier works generally assume bounded subgradients. As was already explained in Part I \cite{ying16ssgd1}, this is a serious limitation (which does not hold even for quadratic risks where the gradient vector is linear in $w$ and grows unbounded). Instead, we shall consider the weaker affine Lipschitz condition (\ref{assume.2}), which was  shown in Part I \cite{ying16ssgd1} to be satisfied automatically by important risk functions such as those arising in popular quadratic,  SVM, and LASSO formulations.

	{\em Notation}: We use lowercase letters to denote vectors, uppercase
	letters for matrices, plain letters for deterministic
	variables, and boldface letters for random variables. We also
	use $(\cdot)^{\sf T}$  to denote transposition, $(\cdot)^{-1}$ for matrix inversion,
	$\mbox{\sf Tr}(\cdot)$ for the trace of a matrix, $\lambda(\cdot)$ for the eigenvalues of
	a matrix, $\|\cdot\|$ for the 2-norm of a matrix or the Euclidean
	norm of a vector, and $\rho(\cdot)$ for the spectral radius of a matrix.
	Besides, we use $ A \geq  B$ to denote that $A - B$ is positive
	semi-definite, and $p\succ 0$ to denote that all entries of vector $p$ are positive.

\section{Problem Formulation: Multi-Agent Case}
	We now extend the single agent scenario analysis to multi-agent networks where a collection of agents cooperate with each other to seek the minimizer of a weighted aggregate cost of the form:
	\be
	\min_{w} \sum_{k=1}^{N}q_k J_k(w) \label{eq.originMulti},
	\ee
	where $k$ refers to the agent index and $q_k$ is some positive weighting coefficient added for generality. When the $\{q_k\}$ are uniform and equal to each other, then (\ref{eq.originMulti}) amounts to minimizing the aggregate sum of the individual risks $\{J_k(w)\}$.
	We can assume, without loss in generality, that the weights $\{q_k\}$ are normalized  to add up to one
	\be
	\sum_{k=1}^N q_k =1   \label{eq.q.norm}
	\ee
	Each individual risk function continues to be expressed as the expected value of some loss function:
	\be
	J_k(w)\define \Ex Q_k(w;\x_k).\ee
	Here, the letter $\x_k$ represents the random data at agent $k$ and the expectation is  over the distribution of this data. Many problems in adaptation and learning involve risk functions of this form, including, for example, mean-square-error designs and support vector machine (SVM) solutions --- see, e.g.,
	\cite{sayed2011adaptive,theodoridis2008Pattern,bishop2006pattern}.  We again allow
	each risk function $J_k(w)$ to be {\em non-differentiable}. This situation is common in machine learning formulations, e.g., in SVM costs and in regularized sparsity-inducing formulations.
	
	{

	We continue to assume that the individual costs satisfy Assumptions~1 and~2 described in the introduction section, namely, conditions (\ref{assump.strongCVX}), (\ref{assume.2}), and (\ref{assump.sbg2}), which ensure that each $J_k(w)$ is strongly-convex and its sub-gradient vectors are affine-Lipschitz with parameters $\{\eta_k,c_k,d_k,e_k,f_k\}$; we are attaching a subscript $k$ to these parameters to make them agent-dependent (alternatively, if desired, we can replace them by agent-independent parameters by using bounds on their values).

}
	\subsection{Network Model}
\begin{figure}[htb]
		\epsfxsize 6.9cm \epsfclipon
		\begin{center}
			\leavevmode \epsffile{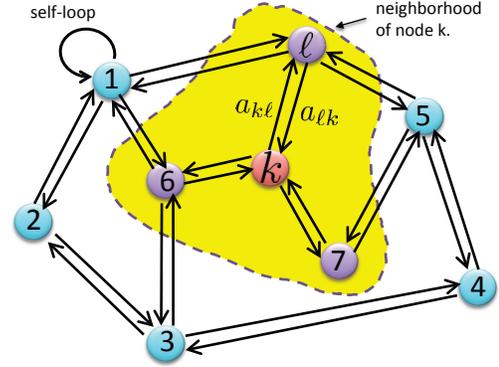} \caption{{\small Agents
					that are linked by edges can share information. The neighborhood of agent
					$k$ is marked by the broken line and consists of the set
					${\cal N}_k=\{6,7,\ell,k\}$.}}\label{fig.network}
		\end{center}\vspace{-3mm}
	\end{figure}

	We consider a network consisting of $N$ separate agents connected by a topology. As described in
	\cite{sayed2014adaptation,sayed2014adaptive}, we assign a pair of nonnegative weights, $\{a_{k\ell},a_{\ell k}\}$, to the edge connecting any two agents $k$ and $\ell$. The scalar $a_{\ell k}$ is used by agent $k$ to scale the data it receives from agent $\ell$ and similarly for $a_{k\ell}$. The network is said to be {\it connected} if paths with nonzero scaling weights can be found linking any two distinct agents in both directions. The network is said to be {\em strongly--connected} if it is connected with at least one self-loop, meaning that $a_{kk}>0$ for some agent $k$.  Figure~\ref{fig.network} shows one example of a strongly--connected network. For emphasis in this figure, each edge between two neighboring agents is represented by two directed arrows. The neighborhood of any agent $k$ is denoted by ${\cal N}_k$ and it consists of all agents that are connected to $k$ by edges; we assume by default that this set includes agent $k$ regardless of whether agent $k$ has a self-loop or not.
	
	\vspace{-1mm}
	There are several strategies that the agents can employ to seek the minimizer, $w^{\star}$, including consensus and diffusion strategies \cite{kar2009distributed,nedic2009distributed,yu2009distributed,sayed2014adaptation,sayed2014adaptive}. As noted earlier, in this work, we focus on the latter class since diffusion implementations have been shown to have superior stability and performance properties over consensus strategies when used in the context of adaptation and learning from streaming data (i.e., when the step-sizes are set to a constant value as opposed to a diminishing value) \cite{sayed2014adaptation,sayed2014adaptive,tu2012diffusion}. Although diminishing step-sizes annihilate the gradient noise term they, nevertheless, disable adaptation and learning in the long run. In comparison, constant step-size updates keep adaptation alive, but they allow gradient noise to seep into the operation of the algorithm. The challenge in these scenarios is therefore  to show that the dynamics of the diffusion strategy over the network is such that the gradient noise effect does not degrade performance and that the network will be able to learn the unknown. This kind of analysis has been answered before in the affirmative for smooth twice-differentiable functions, $J_k(w)$ ---
	see \cite{sayed2014adaptation,sayed2014adaptive,chen2015learning1,chen2015learning2}. In this work, we want to pursue the analysis more generally for possibly {\em non-differentiable} risks in order to encompass important applications (such as SVM learning by multi-agents or LASSO and sparsity-aware learning by similar agents\cite{di2012sparse,di2013sparse,liu2012diffusion,chouvardas2012sparsity}). We also want to pursue the analysis under the weaker affine-Lipschitz assumption (\ref{assume.2}) on the sub-gradients than the stronger conditions used in the prior literature, as we already explained in the earlier sections and in Part I \cite{ying16ssgd1}.
	
	\subsection{Diffusion Strategy}
	We consider the following diffusion strategy in its adapt-then-combine (ATC) form:
	\be
	\left\{\begin{aligned}
		\bpsi_{k,i}&=\w_{k,i-1} - \mu_k\,  \widehat{g}_{k} (\w_{k,i-1}) \\
		\w_{k,i} &=\sum_{\ell\in\mathcal{N}_k} a_{\ell k} \bpsi_{\ell,i}\label{atc.eq2}
	\end{aligned}\right.
	\ee
	Here, the first step involves adaptation by agent $k$ by using a stochastic sub-gradient iteration, while the second step involves aggregation; we assume the gradient noise processes across all agents are independent of each other. The entries $A=[a_{\ell k}]$ define a left-stochastic matrix, namely, the entries of $A$ are non-negative and each of its columns adds up to one. Since the network is strongly-connected, the combination matrix $A$ will be primitive \cite{meyer2000matrix,sayed2014adaptation}. This implies that $A$ will admit a Jordan-decomposition of the form:
	\be
	A = V_{\epsilon} J V_{\epsilon}^{-1} \define \ba{c|c}p & V_R\ea \ba{c|c} 1& 0\\ \hline  0& J_{\epsilon} \ea \ba{c}\one\tran \\ \hline \vspace{-0.3cm}\\ V_L\tran\ea \label{eig.decomp},
	\ee
	with a single eigenvalue at one and all other eigenvalues strictly inside the unit circle. The matrix $J_{\epsilon}$ has a Jordan structure with the ones that would typically appear along its first sub-diagonal replaced by a small positive number, $\epsilon>0$. Note that the eigenvectors of $A$ corresponding to the eigenvalue at one are denoted by
	\be
	Ap=p,\;\;\;A^{\sf T}\one=\one\label{lladk1l3k1l.ad}.
	\ee
	where $\one$ refers to a column vector with all its entries equal to one. It is further known from the Perron-Frobenius theorem \cite{meyer2000matrix} that the entries of $p$ are all strictly positive; we normalize them to add up to one. We denote the individual entries of $p$ by $\{p_k\}$ so that:
	\be
	p_k>0,\;\;\;\;\sum_{k=1}^{N} p_k =1 \label{perron}.
	\ee
	Furthermore, since $ V_{\epsilon}V_{\epsilon}^{-1}= I$, it holds that
	\be
	V_R\tran \one = 0, \quad V_L\tran p = 0,\quad V_L\tran V_R =I \label{two.zero.relation}.
	\ee
	Next, we introduce the vector
	\be q = \mbox{\rm col}\{q_1,q_2,\ldots, q_N \}     \ee
	where $q_k$ is the weight associated with $J_k(w)$ in (\ref{eq.originMulti}).  Since the designer is free to select the step-size parameters, it turns out that we can always relate the vectors $\{p,q\}$ in the following manner:
	\be q = \zeta\; \mbox{\rm diag}\{\mu_1,\mu_2,\ldots,\mu_N\}\;p \label{r3189.f32}\ee
	for some constant  $\zeta >0$. Note, for instance, that for (\ref{r3189.f32}) to be valid the scalar  $\zeta$ should satisfy 
		$\zeta=q_k/\mu_k p_k$ for all $k$. To make this expression for $\zeta$ independent of $k$, we may 
		parameterize (select) the step-sizes as
		\be
		\mu_k = \left(\frac{q_k}{p_k}\right) \mu_o \label{relation}
		\ee
		for some small $\mu_o>0$.  Then, $\zeta=1/\mu_o$, which is independent of $k$ and relation (\ref{r3189.f32}) is satisfied.
	Using (\ref{eq.q.norm}) and (\ref{r3189.f32}) it is easy to check that 
	\be
	\sum_{k=1}^N p_k \mu_k \;=\;\mu_o
	\ee
	Note that since the $\{p_k\}$ are positive, smaller than one, and their sum is one, the above expression shows that $\mu_o$ can be interpreted as a weighted average step-size parameter.

	\section{Network Performance}
	We are now ready to extend Theorem 1 from Part I\cite{ying16ssgd1} to the network case. The analysis is more challenging due to the coupling among the agents. But the result will establish that the distributed strategy is stable and converges exponentially fast for sufficiently small step-sizes. As was the case with Part I \cite{ying16ssgd1}, the  statement below is again in terms of pocket variables, which we define as follows. 
	
	At every iteration $i$, the risk value that is attained by iterate $\w_{k,i}$ is $J_k(\w_{k,i})$. This value is a random variable due to the randomness in the streaming data used to run the algorithm. We denote the mean risk value at agent $k$ by  $\Ex J_k(\w_{k,i})$. We again introduce a {\em best pocket} iterate, denoted by $\w_{k,i}^{\rm best}$. At any iteration $i$, the value that is saved in this pocket variable is the iterate that has generated the smallest mean risk value up to time $i$, i.e.,
	\be
	\w^{\rm best}_{k,i}\;\define\;{\color{black}\argmin_{1\leq j\leq i}}\;\Ex\,J_k(\w_{k,j}).\label{w_k.best}
	\ee
	Observe that in the network case we now have $N$ pocket values, one for each agent. 
	
	\begin{theorem}[{\sc Network performance}]\label{theorem.2}
		Consider using the stochastic sub-gradient diffusion algorithm (\ref{atc.eq2}) to seek the unique minimizer, $w^{\star}$, of the optimization problem (\ref{eq.originMulti}), where the risk functions, $J_k(w)$, are assumed to satisfy assumptions (\ref{assump.strongCVX}), (\ref{assump.sbg2}), and (\ref{usadkh.13lk1l3k}) with parameters $\{\eta_k,\beta_k^2,\sigma_k^2,e_k^2,f_k^2\}$. Assume the step-size parameter is sufficiently small (see condition (\ref{see.dlkarqe})). 
		{\color{black}
		Then, it holds that
		\begin{align}
		 &\hspace{-1mm}\Ex \left(\sum_{k=1}^N q_kJ_k(\w_{k,i}^{\rm best}) - \sum_{k=1}^N q_k J_k(w^\star)\right) \nn\\
		 &\leq  \xi\cdot\alpha^i\sum_{k=1}^Nq_k\Ex\|\w_{k,0}-w^\star\|^2+{}
		 \nn\\
		 &\;\;\hspace{5mm}\frac{\mu_o}{2}\sum_{k=1}^N\Big(q_k f_k^2+q_k\sigma^2_k + 2hq_k\Big[f_k^2+\|g'_k(w^\star)\|^2+\frac{1}{2}\Big]\Big)\nn\\
\label{gj8923.d}
		\end{align}
		The convergence of $\Ex \sum_{k=1}^N q_kJ_k(\w_{k,i}^{\rm best}) $ towards a neighborhood of size $O(\mu_o)$ around $\sum_{k=1}^N q_kJ_k(w^\star)$ occurs at an exponential rate,  $O(\alpha^i)$, dictated by the parameter
		\begin{align}
		\alpha \define& \max_{k}\;\left\{1-\mu_k\big(\eta_k - \mu_o e_k^2  - \mu_o\beta_k^2 - 2\mu_oh e^2_k\big)\right\}\nn\\
		=\;&1-O(\mu_o).
		\end{align}}
		Condition (\ref{see.dlkarqe}) further ahead ensures $\alpha\in(0,1)$. 
	\end{theorem}
	\bp:  The argument is provided in Appendix \ref{app.theorm2}.
	\ep
	
	\smallskip
	
	The above theorem clarifies the performance of the network in terms of the best pocket values across the agents. However, these pocket values are not readily available because the risk values, $J_{k}(\w_{k,i})$, cannot be evaluated. This is due to the fact that the statistical properties of the data are not known beforehand. As was the case with the single-agent scenario in Part I \cite{ying16ssgd1}, a more practical conclusion can be deduced from the statement of the theorem as follows. We again introduce the geometric sum:
	\be
	S_L\define \sum_{j=0}^L \alpha^{L-j}=\alpha S_{L-1}+1\;=\;\frac{1-\alpha^{L+1}}{1-\alpha},
	\ee
	as well as the normalized and convex-combination coefficients:
	\be
	r_L(j) \define\frac{\alpha^{L-j} }{S_L},\;\;\;{\color{black}j=0,1,\ldots,L.} \label{fweofiji23e}
	\ee
	Using these coefficients, we define a weighted iterate at each agent:
	\begin{align}
	\bar\w_{k,L}\define& \sum_{j=0}^L r_L(j)\w_{k,j}\nn\\
	=&\;\;\frac{1}{S_L}\left[\alpha^{L}\w_{k,0}+\alpha^{L-1}\w_{k,1}+\ldots+\w_{k,L}\right] \label{eq.iter.weight}. 
	\end{align}
	and observe that $\bar{\w}_{k,L}$ satisfies the recursive construction:
	\be
	\bar{\w}_{k,L}\;=\;\left(1-\frac{1}{S_L}\right)\bar{\w}_{k,L-1}\;+\;\frac{1}{S_L}\w_{k,L}.
	\ee
	In particular, as $L\rightarrow\infty$, we have $S_L\rightarrow 1/(1-\alpha)$, and the above recursion simplifies in the limit to
	\be
	\bar{\w}_{k,L}\;=\;\alpha \bar{\w}_{k,L-1}\;+\;(1-\alpha)\w_{k,L}.
	\ee
	\begin{corollary}[{\sc Weighted iterates}] \label{corollary.1}Under the same conditions as in Theorem~\ref{theorem.2}, it holds that
		\bq
		&&\hspace{-6mm}\lim_{L\to\infty}\; \Ex \left(\sum_{k=1}^N q_kJ_k(\bar\w_{k,L}) - \sum_{k=1}^N q_k J_k(w^\star)\right)\hspace{8mm}\nn\\
		 &&\leq \frac{\mu_o}{2}\sum_{k=1}^N\Big(q_k f_k^2+q_k\sigma^2_k + 2hq_k\Big[f_k^2+\|g'_k(w^\star)\|^2+\frac{1}{2}\Big]\Big)\nn\\
		 &&=\;O(\mu_o),
		\label{llkad.lk1l3kmjkda}	\eq
		and convergence continues to occur at the same exponential rate, $O(\alpha^L)$.
		
	\end{corollary}
	\bp
	The argument is provided in Appendix \ref{app.corollary1}.
	\ep
	\smallskip

	Result (\ref{llkad.lk1l3kmjkda}) is an interesting conclusion. However, the statement is in terms of the averaged iterate $\bar{\w}_{k,L}$ whose computation requires knowledge of $\alpha$. This latter parameter is a global information, which is not readily available to all agents. Nevertheless, result (\ref{llkad.lk1l3kmjkda}) motivates the following useful distributed implementation with a similar guaranteed performance bound. We can replace $\alpha$ by a design parameter, $\theta$, that is no less than $\alpha$ but still smaller than one, i.e., $\alpha\leq\theta<1$.
	Next, we introduce the  weighted variable:
	\bq
	\bar\w_{k,L}&\define& {\color{black}\sum_{j=0}^L} r_L(j)\w_{k,j},\label{39g.32}
	\eq
	where now \be r_L(j) = \theta^{L-j}/S_L,\;\; {\color{black}j = 0,1\ldots,L,}\ee
	and \be S_L={\color{black}\sum_{j=0}^L\theta^{L-j}.}\ee
	
	\begin{corollary}[{\sc Distributed Weighted iterates}]\label{corollary.2} Under the same conditions as in Theorem~\ref{theorem.2} and $\alpha\leq\theta<1$, relation (\ref{llkad.lk1l3kmjkda}) continues to hold with $\bar{\w}_{k,L}$ in (\ref{eq.iter.weight}) replaced by (\ref{39g.32}). Moreover, convergence now occurs at the exponential rate $O(\theta^L)$.

	\end{corollary}
	%
	\bp
	The argument is similar to the proof of Corollary 2 from Part I \cite{ying16ssgd1}. 
	\ep
	\bigskip
	
	For ease of reference, we summarize in the table below the
	listing of the stochastic subgradient learning algorithm with
	exponential smoothing for which Corollaries \ref{corollary.1} and \ref{corollary.2} hold.\\
	\rule{0.49\textwidth}{1pt}\vspace{-2mm}
	\begin{center} \bf \small Diffusion stochastic subgradient  with exponential smoothing\end{center} \vspace{-4mm}
	\rule{0.49\textwidth}{1pt}
	{\bf \small Initialization}:
	\(
	S_0=1,\; \bar{w}_{k,0}=w_{k,0} = 0,\; \theta=1-O(\mu).  \nn
	\)\\
	{\bf \small repeat for $i\geq 1$}:\\
	\mbox{\bf \small \hspace{6mm}for each agent $k$}:
	\begin{align}
	\bpsi_{k,i}=&\; \w_{k,i-1} -\mu \widehat{g}_k(\w_{k,i-1})\\
	\w_{k,i} =&\; \sum_{\ell\in\mathcal{N}_k}a_{\ell k}\bpsi_{\ell,i}\\
	S_i=&\;\theta S_{i-1} + 1\\
	\bar{\w}_{k,i}=&\;\left(1-\frac{1}{S_i}\right)\bar{\w}_{k,i-1}\;+\;\frac{1}{S_i} \w_{k,i}\hspace{5mm}
	\end{align}
	\mbox{\bf \small\hspace{6mm} end}\\
	{\bf \small end}\\
	\rule{0.49\textwidth}{1pt}
	
	\smallskip
	\subsection{Interpretation of Results}
	\noindent Examining the bound in (\ref{llkad.lk1l3kmjkda}), and comparing it with result (88) from Part I\cite{ying16ssgd1} for the single-agent case, we observe that the topology of the network is now reflected in the bound through the weighting factor, $q_k$ and step-size $\mu_k$, which can be related to the Perron entry $p_k$ through (\ref{relation}). Recall from (\ref{lladk1l3k1l.ad}) that the $\{p_k\}$ are the entries of the right-eigenvector of $A$ corresponding to the eigenvalue at one. Moreover, the bound in (\ref{llkad.lk1l3kmjkda})involves three terms (rather than only two as in the single-agent case --- compared with (88) from Part I\cite{ying16ssgd1}):
	\begin{enumerate}
		\item[(1)] $q_k f_k^2$, which  arises from the non-smoothness of the risk function;
		\item[(2)] $q_k\sigma^2_k$, which is due to gradient noise and the approximation of the true sub-gradient vector;
		\item[(3)] $2q_kh\Big[f_k^2+\|g'_k(w^\star)\|^2+\frac{1}{2}\Big]$, which is an extra term in comparison to the single agent case. We explained in (\ref{kjahd6713.aldklakd}) that the value of $h$ is related to how far the error at each agent is away from the weighted average error across the network. As for $\|g'_k(w^\star)\|^2$, this quantity represents the disagreement among the agents over $w^\star$. Because each function $J_k(\cdot)$ may have a different minimizer, $g'_k(w^\star)$ is generally nonzero.
	\end{enumerate}
	
	\section{Simulations}
\begin{figure}[ht]
		\centering
		\includegraphics[scale=0.37]{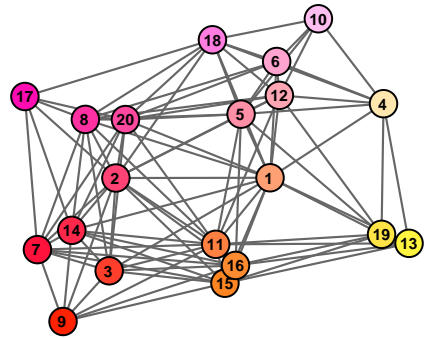}
		\caption{{Network topology linking $N=20$ agents.}} 
		\label{fig.topology}
	\end{figure}

\begin{figure}[ht]
		\centering
		\includegraphics[scale=0.37]{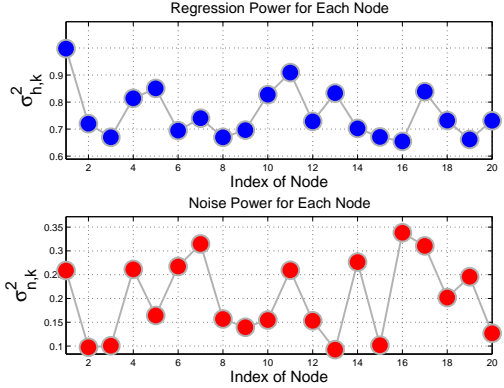} \caption{{Feature and noise variances across the agents. }} \label{fig.SNR}
	\end{figure}

\begin{figure}[ht]
		\epsfxsize 8.1cm \epsfclipon
		\begin{center}
			\leavevmode \epsffile{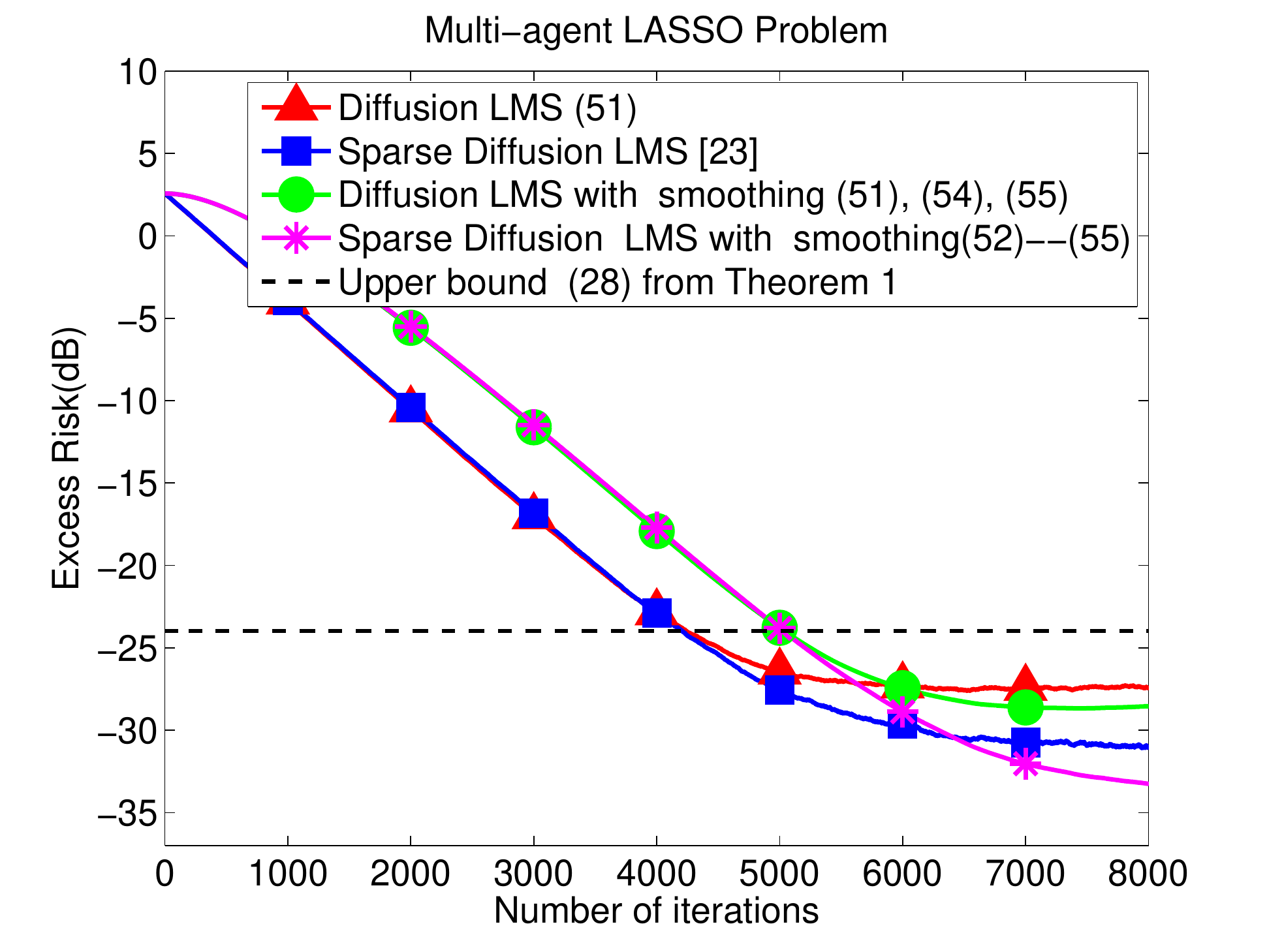} \caption{{\small The excess-risk curves for several strategies.}}\label{fig.LASSO_Network}
		\end{center}
	\end{figure}

	\noindent {\bf Example 1 (Multi-agent LASSO problem)} We now consider the LASSO problem with 20 agents connected according to Fig.~\ref{fig.topology}. A quick review of the LASSO problem is as follows. (A more detailed discussion and the relationship between the proposed assumptions \eqref{assump.strongCVX}--\eqref{assump.sbg2} and the LASSO formulation can be found in Part I \cite{ying16ssgd1}.) We consider follwing cost function for each agent:
	\be
		J_k^{\rm lasso}(w) \define \frac{1}{2}\Ex\|\bgamma_k - \h_k\tran w \|^2 +\delta\|w\|_1,\label{lkadh,13k1jl3k}
	\ee
	where $\delta>0$ is a regularization parameter and $\|w\|_1$ denotes the $\ell_1-$norm of $w$. The variable $\bm{\gamma}_k$ plays the role of a desired signal for agent $k$, while $\h_k$ plays the role of a regression vector for the same agent.  It is assumed that the regression data are zero-mean wide-sense stationary, and its distribution satisfies the standard Gaussian distribution, i.e., $\h_k\sim \mathcal{N}(0,\sigma_{h,k}^2I)$. We further assume that $\{\bm{\gamma}_k,\h_k\}$ satisfy a linear model of the form $\bgamma_k$ generated through:
	\be
		\bgamma_k = \h_k\tran w^o_k +\n_k 
	\ee
	where $\n_k \sim \mathcal{N}(0,\sigma_{n,k}^2I)$ and $w^o_k$ is some sparse random model for each agent.
	Each agent is allowed to have different regression and noise powers, as illustrated in Fig.~\ref{fig.SNR}. 
	Under these modeling assumptions, we can determine a closed-form expression for $w^{\star}$ as follows:
	\begin{align}
		w^\star\hspace{-1mm}
		=& \argmin_{w} \sum_{k=1}^N q_k J_k(w)\nn\\
		=&\, \argmin_{w}\frac{1}{2}\sum_{k=1}^N q_k\sigma_{h,k}^2\|w-w_k^o\|^2  + \delta\|w\|_1\nn\\ 
		=&\,\argmin_{w}\frac{1}{2}\sum_{k=1}^N q_k\sigma_{h,k}^2 \|w\|^2 -\sum_{k=1}^N q_k\sigma_{h,k}^2[w_k^o]\tran w  + \delta\|w\|_1
		\nn\\
		=&\, \argmin_{w}\frac{1}{2}\hspace{-1mm}\left(\sum_{k=1}^N q_k \sigma_{h,k}^2\right)\hspace{-1.mm}\left\|w-\frac{\sum_{k=1}^N q_k\sigma_{h,k}^2 w^o_k}{\sum_{k=1}^N q_k\sigma_{h,k}^2}\right\|^2 \hspace{-1.5mm}+\hspace{-0.5mm} \delta\|w\|_1
	\end{align}
	From first-order optimality conditions, we obtain\cite{donoho1994ideal}:
	\be
		w^\star = \mathcal{S}_\epsilon\left(\frac{\sum_{k=1}^Nq_k \sigma_{h,k}^2 w^o_k}{\sum_{k=1}^Nq_k \sigma_{h,k}^2}\right),
	\ee
	where the symbol $\mathcal{S}_\epsilon$ represents the soft-thresholding function with parameter $\epsilon$, i.e.,
	\be
	\mathcal{S}_\epsilon(x) ={\rm sgn}(x)\cdot \max\{0,|x| - \epsilon\}.
	\ee
	and 
	\be
		\epsilon = \frac{\delta}{\sum_{k=1}^Nq_k \sigma_{h,k}^2} 
	\ee
	where the notation ${\rm sgn}(a)$, for a scalar $a$, refers to the sign function:
	\be
	\mbox{\rm sgn}[a]=\left\{\begin{array}{ll}+1,\;\;&a > 0\\
		\;\;\,0, \;\;&a = 0\\
		-1,&a<0\end{array}\right.
	\ee
	For the stochastic sub-gradient implementation, the following instantaneous approximation for the sub-gradient is employed:
	\begin{align}
		\widehat{g}_k^{\rm lasso}(\w_{i-1}) &= -\h_{k,i}(\bgamma_k(i) - \h_{k,i}\tran \w_{k,i-1}) + \delta\cdot {\rm sgn}(\w_{k, i-1})
		\label{isthel;alkdak}
	\end{align}
	 In Fig. \ref{fig.LASSO_Network}, we compare the performance of this solution against several strategies including standard diffusion LMS\cite{cattivelli2010diffusion,sayed2014adaptation,sayed2014adaptive}:
	\be
	\left\{
	\begin{aligned}
		\bpsi_{k,i}=&\w_{k,i-1} + \mu \h_{k,i}(\bgamma_k(i)-\h_{k,i}\tran\w_{k,i-1})\\
		\w_{k,i} =&\sum_{\ell\in\mathcal{N}_k} a_{\ell k} \bpsi_{\ell,i}
	\end{aligned}
	\right.\label{diffusion}
	\ee
	and sparse diffusion  LMS\cite{di2012sparse,liu2012diffusion,chouvardas2012sparsity} \cite[Eq. 21]{di2013sparse}.\\
	\rule{0.49\textwidth}{1pt}\vspace{-2mm}
	\begin{center} \bf \small Diffusion sparse LMS with expoential smoothing \end{center} \vspace{-4mm}
	\rule{0.49\textwidth}{1pt}
	{\bf \small Initialization}:
	\(
	S_0=1,\; \bar{\w}_{k,0}=\w_{k,0} = 0,\; \theta=1-O(\mu).  \nn
	\)\\
	{\bf \small repeat for $i\geq 1$}:\\
	\mbox{\bf \small \hspace{6mm}for each agent $k$}:
	\begin{align}
	\bpsi_{k,i}=&\; \w_{k,i-1} + \mu_k \h_{k,i}(\bgamma_k(i)-\h_{k,i}\tran\w_{k,i-1}) \nn\\
	&\hspace{3mm}{}-\mu_k\delta\cdot{\rm sgn}(\w_{k,i-1})\\
	\w_{k,i} =&\; \sum_{\ell\in\mathcal{N}_k}a_{\ell k}\bpsi_{\ell,i} \label{238j.gj8}
	\\
	S_i=&\;\theta S_{i-1} + 1 \label{8gm.2j9}\\
	\bar{\w}_{k,i}=&\;\left(1-\frac{1}{S_i}\right)\bar{\w}_{k,i-1}\;+\;\frac{1}{S_i} \w_{k,i}\hspace{5mm} \label{238jg.9}
	\end{align}
	\mbox{\bf \small\hspace{6mm} end}\\
	{\bf \small end}\\
	\rule{0.49\textwidth}{1pt}

	The parameter setting is as follows: $w_k^o\in \real^{100}$ has 5 random non-zero entries uniformly distributed between 0.5 and 1.5, and $\delta=0.005$. We simply let $q_k=p_k$ and set the step-size for all agents at $\mu_k=\mu_o = 0.001$. From the simulations we find $h=1.24$ for the factor that appears in (\ref{gj8923.d}). As for the exponential smoothing factor $\theta$, we chose $\theta = 1 - 2\mu_o(\frac{1}{N}\sum_{k=1}^N\eta_k)=0.9985$.
	
	$\hfill \Box$
	
	
	\noindent {\bf Example 2 (Multi-agent SVM learning)}
	Next, we will consider the multi-agent SVM problem. Similar to LASSO problem, we provide a brief review for notation. More detailed discussion can be found in Part I \cite{ying16ssgd1}.
	The regularized SVM risk function for each agent is of the form:
	\be
		J_k^{\rm svm}(w) \define \frac{\rho}{2} \| w\|^2 + \Ex \left( \max \left\{0, 1-\bgamma_k  \h_k\tran w \right\} \right),
	\ee
	where $\rho>0$ is a regularization parameter. We are generally given a collection of independent training data, $\{\bm{\gamma}_k(i),\h_{k,i}\}$, consisting of feature vectors and their class designations. We select $q_k=\frac{1}{N}$ and
	\be
		\mu_k = \mu_o/(Np_k)
	\ee
	One approximation for the sub-gradient construction at a generic location $\w$ corresponding to generic data $\{\bm{\gamma},\h\}$ is
	\be
	\widehat{g}^{\rm svm}(\w) = \rho \w + \bgamma \h\, \mathds{I} [\bgamma\h\tran \w \leq 1],  \label{BBB}
	\ee
	where the indicator function $\mathbb{I}[a]$ is defined as follows:
	\be
	\mathds{I}[a]\;=\;\left\{\begin{array}{ll}1,&\mbox{\rm if statement $a$ is true}\\0,&\mbox{\rm otherwise}\end{array}\right.
	\ee\newpage
	\noindent\rule{0.49\textwidth}{1pt}\vspace{-2mm}
	\begin{center} \bf \small Diffusion SVM with exponential smoothing\end{center} \vspace{-4mm}
	\rule{0.49\textwidth}{1pt}
	{\bf \small Initialization}:
	\(
	S_0=1,\; \bar{\w}_{k,0}=\w_{k,0}=0,\; \theta=1-O(\mu).  \nn
	\)\\
	{\bf \small repeat for $i\geq 1$}:\\
	\mbox{\bf \small \hspace{6mm}for each agent $k$}:
	\begin{align}
	\bpsi_{k,i}=&\; (1-\rho\mu)\w_{k,i-1} \nn\\
	&\;\;\;{}- \mu\bgamma_k(i)\h_i\mathbb{I}[\bgamma(k,i)\h_{k,i}\tran\w_{k,i-1}\leq 1] \\
	\w_{k,i} =&\; \sum_{\ell\in\mathcal{N}_k}a_{\ell k}\bpsi_{\ell,i}\\
	S_i=&\;\theta S_{i-1} + 1\\
	\bar{\w}_{k,i}=&\;\left(1-\frac{1}{S_i}\right)\bar{\w}_{k,i-1}\;+\;\frac{1}{S_i} \w_{k,i}\hspace{5mm}
	\end{align}
	\mbox{\bf \small\hspace{6mm} end}\\
	{\bf \small end}\\
	\rule{0.49\textwidth}{1pt}\vspace{1mm}
	
\begin{figure}[ht]
		\centering
			\includegraphics[scale=0.37]{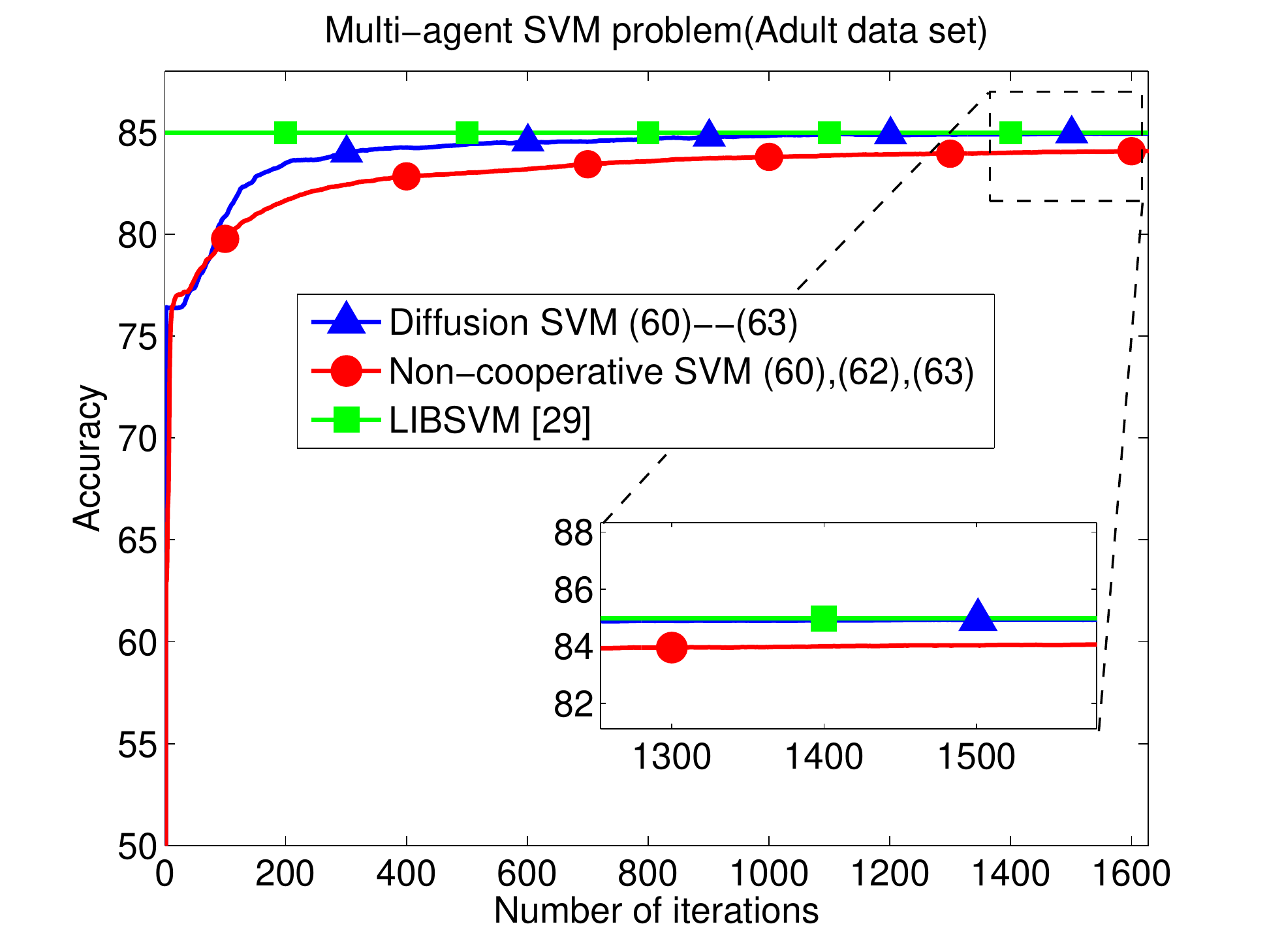} 
			\includegraphics[scale=0.37]{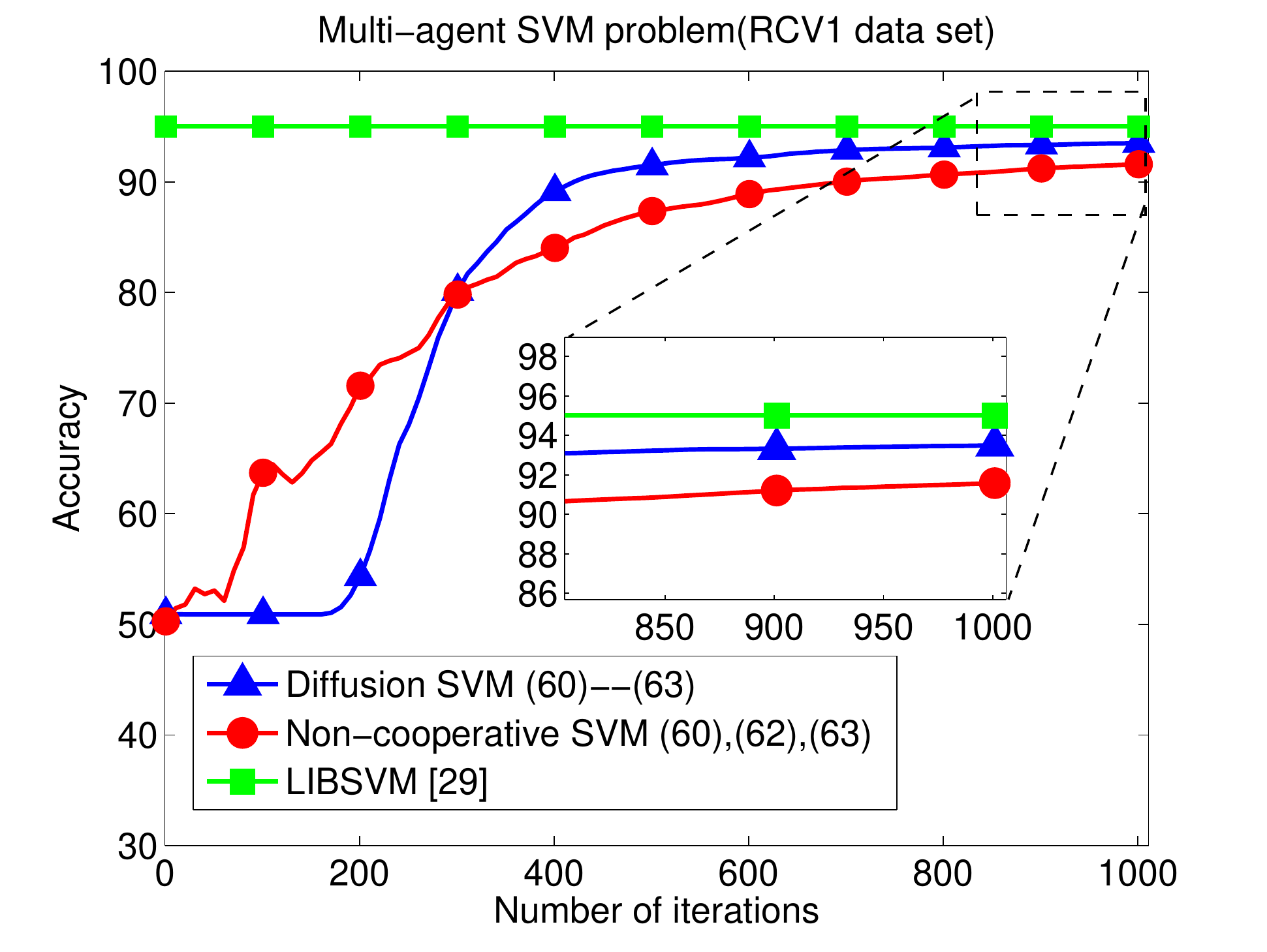}
			\caption{{\small Performance of diffusion SVM for the Adult dataset (Top) and RCV1 dataset (Bottom), where vertical axis measures the percentage of correct prediction over test dataset.}} \label{fig.SVM_Adult_Network}
	\end{figure}

	We distribute 32561 training data from an adult dataset\footnote{\url{https://archive.ics.uci.edu/ml/datasets/Adult}} over a network consisting of 20 agents. We set $\rho=0.002$ and $\mu_o=0.15$ for all agents. From Example 6 in Part I \cite{ying16ssgd1} and Theorem \ref{theorem.2}, we know that for the multi-agent SVM problem:
	\bq
	\alpha &=& \max_k\; \left\{1-\mu\rho+\mu^2(2h+1)e_k^2\right\}\nn\\
	&=& \max_k\; \left\{1-\mu\rho+\mu^2(2h+1)2\rho^2\right\}.
	\eq
	We set $\theta = 1-0.9\cdot\mu_o\rho$, which usually guarantees $\theta\geq \alpha$. Fig. \ref{fig.SVM_Adult_Network} (left) shows that cooperation among the agents outperforms the non-cooperative solution. Moreover, the  distributed network can almost match the performance of the centralized LIBSVM solution\cite{CC01a}.  We also examined the RCV1 dataset\footnote{\url{https://www.csie.ntu.edu.tw/~cjlin/libsvmtools/datasets/binary.html}}. Here we have 20242 training data points and we distribute them over 20 agents. We set the parameters to $\rho=1\times 10^{-5}$ and $\mu_o=0.5$ (due to limited data). We now use $\theta = 1-0.5\cdot\mu_o\rho$ since $\mu$ is not that small. The result is shown in Fig.~\ref{fig.SVM_Adult_Network} (right).
	
	$\hfill \Box$


	\section{Conclusion}
	In summary, we examined the performance of stochastic sub-gradient learning strategies over adaptive networks. We proposed a new affine-Lipschitz condition, which is quite suitable for strongly convex but non-differentiable cost functions and is automatically satisfied by   several important cases including SVM, LASSO, Total-Variation denoising, etc. Under this weaker condition, the analysis establishes that sub-gradient strategies can attain exponential convergence rates, as opposed to sub-linear rates. The analysis also establishes that these strategies can approach the optimal solution within $O(\mu)$, for sufficiently small step-sizes.
	\\

	\appendices
	\section{\sc Proof of theorem \ref{theorem.2}}\label{app.theorm2}
	Introduce the error vector, $\widetilde{\w}_{k,i}=w^{\star}-\w_{k,i}$. We collect the iterates and the respective errors from across the network into block column vectors:
	\bq	\swb_{i} &\define &\; {\rm col} \{\w_{1,i},\w_{2,i},\hdots,\w_{N,i} \}\\
	\widetilde{\swb}_{i}&\define &\; {\rm col} \{\widetilde{\w}_{1,i},\widetilde{\w}_{2,i},\hdots,\widetilde{\w}_{N,i}\}.
	\eq
	We also define the extended quantities:
	\begin{align}
		\mathcal{A} &\define \; A\otimes I_M\\
		\sgb(\swb_{i-1} ) &\define\; {\rm col}\{g_1(\w_{1,i-1}),,\hdots,g_N(\w_{N,i-1})\} \\
		\ssb_{i}(\swb_{i-1}) &\define \;{\rm col}\{s_{1,i}(\w_{1,i-1}),\hdots,s_{N,i}(\w_{N,i-1})\},\\
		U &\define\; {\rm diag}\left\{\mu_1,\mu_2\cdots, \mu_N\right\} / \mu_o\\
		\mathcal{U}  &\define\; U\otimes I_M
	\end{align}
	where $\otimes$ denotes the Kronecker product operation, and $\s_{k,i}(\w_{k,i-1})$ denotes the gradient noise at agent $k$.
	Using this notation, it is straightforward to verify that the network error vector generated by the diffusion strategy \eqref{atc.eq2} evolves according to the following dynamics:
	\be
	\widetilde\swb_{i} = \mathcal{A}\tran \left(\widetilde\swb_{i-1} + \mu_o \mathcal{U}\sgb(\swb_{i-1} ) +\mu_o \mathcal{U} \ssb_{i}(\swb_{i-1}) \right)\label{eq.net.errorRec1}.
	\ee
	Motivated by the treatment of the smooth case in \cite{sayed2014adaptation,chen2015learning1,chen2015learning2}, we introduce a useful change of variables. Let ${\cal V}_{\epsilon}=V_{\epsilon}\otimes I_M$ and
	${\cal J}_{\epsilon}=J_{\epsilon}\otimes I_M$. Multiplying (\ref{eq.net.errorRec1}) from the left by ${\cal V}_{\epsilon}^{\sf T}$ gives
	\begin{align}
	\mathcal{V}_{\epsilon}\tran\widetilde\swb_{i} =&\; \mathcal{J}\tran\left[\mathcal{V}_{\epsilon}\tran\widetilde\swb_{i-1} +
	\mu_o\mathcal{V}_{\epsilon}\tran \cU \sgb(\swb_{i-1})\right. \nn\\
	&\;\;\;\;\;\;\;\;\left.{}+\mu_o \mathcal{V}_{\epsilon}\tran \cU\ssb_{i}(\swb_{i-1})\right].
	\label{refalkdl.13lk}\end{align}
	where from \eqref{eig.decomp}:
	\be
		\cJ \define  \ba{c|c} 1& 0\\ \hline  0& J_{\epsilon} \ea  \otimes I_M
	\ee
	and 
	\begin{align}
		\mathcal{V}_{\epsilon}\tran \cU =&\, \left(\ba{c}p\tran\\ V_R\tran\ea\otimes I_M \right)(U\otimes I_M)\nn\\
		=&\, \ba{c}p\tran U \\ V_R\tran U\ea\otimes I_M \nn\\
		\stackrel{(\ref{r3189.f32})}{=}&\, \ba{c} q\tran\otimes I_M \\ V_R\tran U \otimes I_M\ea
	\end{align}
	To proceed, we introduce 
	\begin{align}
		\mathcal{V}_{\epsilon}\tran\widetilde\swb_{i} \hspace{-1mm}\,=\,\hspace{-1mm} \ba{c} \hspace{-1mm}(p\tran \otimes I)\widetilde\swb_{i}\hspace{-1mm}\\ (V_R\tran \otimes I) \widetilde \swb_{i} \ea & \hspace{-1mm}\define\hspace{-1mm} \ba{c} \bar\w_{i} \\ \check\swb_{i} \ea,\nn\\\label{def.w2}
		\\
		\mathcal{V}_{\epsilon}\tran\cU  \sgb(\swb_{i-1})\hspace{-1mm}=\hspace{-1mm} \ba{c}\hspace{-1mm} (q\tran \otimes I)\sgb(\swb_{i-1}) \hspace{-3mm}\\\hspace{-1mm} (V_R\tran U \otimes I) \sgb( \swb_{i-1}) \hspace{-3mm}\ea& \hspace{-1mm}\define\hspace{-1mm} \ba{c}\bar{g}(\swb_{i-1})\\ \check{\sgb}(\swb_{i-1}) \ea\label{def.transSG}
		\\
		\mathcal{V}_{\epsilon}\tran\cU \ssb_{i}(\swb_{i-1})\hspace{-1mm}=\hspace{-1mm}\ba{c} (q\tran \otimes I)\ssb_{i}(\swb_{i-1})\\\hspace{-1mm} (V_R\tran U \otimes I) \ssb_{i}( \swb_{i-1})\hspace{-1mm} \ea& \hspace{-1mm}\define\hspace{-1mm} \ba{c}\bar \s_{i}(\swb_{i-1}) \\ \check \ssb_{i}(\swb_{i-1}) \ea \label{def.transNoise}
	\end{align}
	where the quantities $\{\bar\w_{i}, \bar{g}(\swb_{i-1}), \bar \s_{i}(\swb_{i-1})\}$ amount to the weighted averages:
	\bq
	\bar\w_{i}&=&\sum_{k=1}^N p_k \widetilde{\w}_{k,i},\\
	\bar{g}(\swb_{i-1})&=&\sum_{k=1}^N q_k g_k(\w_{k,i-1}),\\
	\bar \s_{i}(\swb_{i-1})&=&\sum_{k=1}^N q_k \s_{k,i}(\w_{k,i-1}).
	\eq
	It is useful to observe the asymmetry reflected in the fact that $\bar{\w}_i$ is obtained by using the weights $\{p_k\}$ while the averages \eqref{def.transSG}--\eqref{def.transNoise} are obtained by using the weights $\{q_k\}$. 
	We can now rewrite (\ref{refalkdl.13lk}) as
	\begin{align}
		\ba{c} \bar\w_{i} \\ \check\swb_{i} \ea=&\;\ba{c|c}  I_M& 0\\ \hline  0& \mathcal{J}\tran_{\epsilon} \ea \left(\ba{c} \bar\w_{i-1} \\ \check\swb_{i-1} \ea \right.\label{eq.errorReccursion2}\\
		& \;\;\left.{}+ \mu_o\ba{c}\bar{g} (\swb_{i-1})\\ \check{\sgb}(\swb_{i-1}) \ea\right.
		\left.+\mu_o\ba{c}\bar{\s}_{i}(\swb_{i-1}) \\ \check{\ssb}_{i}(\swb_{i-1}) \ea \right).\nn 
	\end{align}
	Consider the top recursion, namely,
	\be
	\bar\w_{i} =  \bar\w_{i-1} + \mu_o \bar{g}(\swb_{i-1}) + \mu_o \bar{\s}_{i}(\swb_{i-1}). \label{eq.barErrRec}
	\ee
	Squaring and taking expectations we have
	\begin{align}
		&\hspace{-9mm}\Ex [\|\bar\w_{i}\|^2 \,|\,\bm{\cal F}_{i-1}]\nn\\
		\hspace{2mm}=&\; \Ex [\|\bar\w_{i-1} + \mu_o \bar{g}(\swb_{i-1}) +\mu_o\bar{\s}_{i}(\swb_{i-1})\|^2  \,|\,\bm{\cal F}_{i-1}] \nn\\
		=&\; \|\bar\w_{i-1} \|^2 + 2\mu_o\bar g(\swb_{i-1})\tran \bar\w_{i-1}\;+\mu_o^2  \| \bar g(\swb_{i-1})\|^2 \;\nn\\
		&\;\;{}+\; \mu_o^2\Ex [\|\bar \s_{i}(\swb_{i-1})\|^2 \,|\,\bm{\cal F}_{i-1}] \label{eq.BarErr1}.
	\end{align}
	We examine the terms on the right-hand side one by one. First note that, using Jensen's inequality,
	{\color{black}
	\begin{align}
	\| \bar g(\swb_{i-1})\|^2 =&\;  \left\|\sum_{k=1}^Nq_k g_k(\w_{k,i-1})\right\|^2 \nn\\
	\stackrel{(a)}{=}&\; \left\| \sum_{k=1}^N q_k g_k(\w_{k,i-1}) - \sum_{k=1}^N q_k g'_k(w^\star) \right\|^2\nn\\
	\leq&\; {\color{black}\sum_{k=1}^Nq_k \|g_k(\w_{k,i-1}) - g'_k(w^\star)\|^2 }\nn\\
	\stackrel{(\ref{assump.sbg2})}{\leq}&\; \sum_{k=1}^Nq_k \Big( e_k^2 \| \widetilde \w_{k,i-1} \|^2 + f_k^2 \Big) \label{eq.ineq1}.
	\end{align}
	In step (a), we exploit the fact that, by definition, $w^{\star}$ is the minimizer of (\ref{eq.originMulti}) and, hence, there exist sub-gradients $g'_k(w^\star), k=1,2,\cdots,N$, satisfying $\sum_{k=1}^N q_kg'_k(w^\star)=0$. }
	
	Next, the noise term can be bounded by:
	\begin{align}
		\Ex[ \| \bar{\s}_{i}(\swb_{i-1})\|^2  \,|\,\bm{\cal F}_{i-1}] &= \Ex \left[ \left\|\sum_{k=1}^Nq_k \s_k(\w_{k,i-1})\right\|^2 \,|\,\bm{\cal F}_{i-1}\right]  \nn\\
		&\stackrel{(a)}{\leq} \sum_{k=1}^Nq_k\Ex [\| \s_k(\w_{k,i-1})\|^2  \,|\,\bm{\cal F}_{i-1}] \nn\\
		&\leq \sum_{k=1}^Nq_k\Big( \beta_k^2 \| \widetilde \w_{k,i-1} \|^2 + \sigma_k^2\Big) \label{eq.ineq2}.
	\end{align}
	where step (a) follows from  Jensen's inequality. 
	
	Finally, with regards to the cross term in (\ref{eq.BarErr1}), we adapt  an argument from \cite{nedic2009distributed} to obtain (\ref{eq.disgust1}) by first noting that:
	\begin{align}
		&\hspace{-2mm}\bar g(\swb_{i-1})\tran \bar\w_{i-1} \nn\\
		&= \sum_{k=1}^Nq_k g_k\tran(\w_{k,i-1}) \big(\widetilde \w_{k,i-1}+\bar\w_{i-1}-\widetilde \w_{k,i-1}\big) \nn\\
		&=\sum_{k=1}^Nq_k g_k\tran(\w_{k,i-1}) \widetilde \w_{k,i-1}
		\nn\\
		&\hspace{6mm}+\sum_{k=1}^Nq_k g_k\tran(\w_{k,i-1}) \big(\bar\w_{i-1}-\widetilde \w_{k,i-1}\big).
		\label{3h89g.8hg}
	\end{align}
	Using the strong-convexity property (\ref{assump.strongCVX}), we have 
	\be
	g_k(\w_{k,i-1})\tran \widetilde{\w}_{k,i-1} \leq J_k(w^\star) - J_k(\w_{k,i-1}) - \frac{\eta_k}{2} \|\widetilde{\w}_{k,i-1}\|^2,
	\ee
	Substituting into (\ref{3h89g.8hg}) gives
	\bq 
	&&\hspace{-1cm}\bar g(\swb_{i-1})\tran \bar\w_{i-1} 	\nn\\
	&\leq&\sum_{k=1}^Nq_k \Big(J_k(w^\star) - J_k(\w_{k,i-1}) - \frac{\eta_k}{2} \|\widetilde{\w}_{k,i-1}\|^2\Big) +\nn\\
	&&\sum_{k=1}^Nq_k g_k\tran(\w_{k,i-1})\big(\bar\w_{i-1}-\widetilde \w_{k,i-1}\big) \nn \\
	&\leq&\sum_{k=1}^Nq_k \Big(J_k(w^\star) - J_k(\w_{k,i-1}) - \frac{\eta_k}{2} \|\widetilde{\w}_{k,i-1}\|^2\Big) +\nn\\
	&&\sum_{k=1}^Nq_k \|g_k(\w_{k,i-1})\|\|\bar\w_{i-1}-\widetilde \w_{k,i-1}\|.\label{eq.disgust1}
	\eq
	It follows, under expectation, that
	\bq  
	&&\hspace{-1cm}\Ex\bar g(\swb_{i-1})\tran \bar\w_{i-1}\nn\\
	&\leq&\sum_{k=1}^Nq_k \Big(J_k(w^\star) - \Ex J_k(\w_{k,i-1}) - \frac{\eta_k}{2} \Ex\|\widetilde{\w}_{k,i-1}\|^2\Big) +\nn\\
	&&\sum_{k=1}^Nq_k\Ex\left( \|g_k(\w_{k,i-1})\|\|\bar\w_{i-1}-\widetilde \w_{k,i-1}\|\right).
	\eq
	Now, using the Cauchy-Schwartz inequality, we can bound the last expectation as
	\begin{align}
		\Ex\Big(\|g_k(\w_{k,i-1})\| \|\bar\w_{i-1}-\widetilde \w_{k,i-1}\|\Big)\hspace{2cm} \nn\\
		\leq\sqrt{\Ex \|g_k(\w_{k,i-1})\|^2 \Ex \|\bar\w_{i-1}-\widetilde \w_{k,i-1}\|^2}.
	\end{align}
	After sufficient iterations, it will hold that (see  Appendix \ref{app.netClose} for the proof):
	\be
	\Ex\|\bar\w_{i-1}-\widetilde \w_{k,i-1}\|^2 = O(\mu_o^2)\label{prove.this.appendix}.
	\ee
	This means that there exists an $I_o$ large enough and a constant $h$  such that for all $i\geq I_o$:
	\be
	\Ex\|\bar\w_{i-1}-\widetilde \w_{k,i-1}\|^2 \leq h^2\mu_o^2.
	\label{kjahd6713.aldklakd}\ee
	Therefore, we find that
	\begin{align}
		&\hspace{-0.8cm}\Ex\Big(\|g_k(\w_{k,i-1})\|\|\bar\w_{i-1}-\widetilde \w_{k,i-1}\|\Big) \nn\\
		\leq&\;h\mu_o\,\left(\sqrt{\Ex \|g_k(\w_{k,i-1})\|^2}\right)  \nn\\
		\leq&\; h\mu_o\,\left(\sqrt{2\Ex \|g_k(\w_{k,i-1})-g_k'(w^\star)\|^2 + 2\|g_k'(w^\star)\|^2}\right)  \nn\\
		\stackrel{(\ref{assump.sbg2})}{\leq}\hspace{-1mm}&\; h\mu_o\;\left(\sqrt{2e_k^2\Ex\|\widetilde\w_{k,i-1}\|^2+2f_k^2+2\|g_k'(w^\star)\|^2}\right) \nn\\
		\leq&\; h\mu_o\left(\frac{e_k^2\Ex\|\widetilde\w_{k,i-1}\|^2+f_k^2+\|g_k'(w^\star)\|^2}{R}+\frac{R}{2} \right) \label{eq.128},
	\end{align}
	where the last inequality follows from using
	\be
	\sqrt{x} \leq \frac{1}{2}\Big(\frac{x}{R}+R\Big),\;\;x\geq 0 \label{sqrt.ineq},
	\ee
	which follows from the inequality
		\be
		 \frac{1}{2}\frac{x}{R}-\sqrt{x} + \frac{1}{2}R=\frac{1}{2}\left(\sqrt{\frac{x}{R}}-\sqrt{R}\right)^2\geq0 \nn
		\ee
	for any positive $R$, e.g., $R=1$, which allows us to conclude that, as $i\rightarrow\infty$:
	\begin{align}
		&\hspace{-0.6cm}\Ex\bar g(\swb_{i-1})\tran \bar\w_{i-1} \nn\\
		\leq&\; \sum_{k=1}^Nq_k \Big(J_k(w^\star) - \Ex J_k(\w_{k,i-1}) - \frac{\eta_k}{2} \Ex\|\widetilde{\w}_{k,i-1}\|^2\Big) \nn\\
		&{}+ \mu_o  \sum_{k=1}^Nhq_k \left( {e_k^2\Ex\|\widetilde\w_{k,i-1}\|^2+f_k^2+\|g_k'(w^\star)\|^2}+\frac{1}{2} \right)   
		\nn\\
		 \label{eq.disgust2}
	\end{align}
	Taking expectation of (\ref{eq.BarErr1}) over the filtration and substituting (\ref{eq.ineq1}), (\ref{eq.ineq2}), and (\ref{eq.disgust2}), we obtain asymptotically that:
	\bq&&\hspace{-1cm} \Ex \|\bar\w_i\|^2\nn\\
	 \hspace{-2mm}&\leq &\hspace{-2mm}\Ex\|\bar\w_{i-1}\|^2 +2\mu_o\sum_{k=1}^Nq_k \Big(J_k(w^\star) - \Ex J_k(\w_{k,i-1}) \Big) \nn\\
	 \hspace{-2mm}&&\hspace{-1mm}{}- \mu_o\sum_{k=1}^Nq_k\eta_k \Ex\|\widetilde{\w}_{k,i-1}\|^2 \nn\\
	 \hspace{-2mm}&& \hspace{-1mm}
	 {}+\mu_o^2\sum_{k=1}^Nq_k \Big( e_k^2\Ex \| \widetilde \w_{k,i-1} \|^2 + f_k^2 \Big)\nn\\
	 \hspace{-2mm}&&\hspace{-1mm}{}+\mu_o^2\sum_{k=1}^Nq_k\Big( \beta_k^2\Ex \| \widetilde \w_{k,i-1} \|^2 + \sigma_k^2\Big) \nn\\
	 \hspace{-2mm}&&\hspace{-1mm}{}+ 2\mu_o^2 \sum_{k=1}^Nq_kh\hspace{-1mm}\left(\hspace{-0.5mm} {e_k^2\Ex\|\widetilde\w_{k,i-1}\|^2+f_k^2+\|g_k'(w^\star)\|^2}\hspace{-0.2mm}+\hspace{-0.2mm}\frac{1}{2} \hspace{-0.5mm}\right)   \nn\\
	 \hspace{-2mm}&\leq&\hspace{-2mm}\Ex\|\bar\w_{i-1}\|^2 + 2\mu_o\sum_{k=1}^Nq_k \Big(J_k(w^\star) - \Ex J_k(\widetilde\w_{k,i-1}) \Big) \nn\\
	\hspace{-2mm} & &\hspace{-1mm}{}-\sum_{k=1}^N\left(1-\alpha_k\right) p_k\Ex\|\widetilde \w_{k,i-1}  \|^2\nn\\
	 \hspace{-2mm}&&\hspace{-1mm} {}+\mu_o^2 \sum_{k=1}^N\Big(q_k f_k^2+q_k\sigma^2_k + 2hq_k\Big[f_k^2+\|g'_k(w^\star)\|^2+\frac{1}{2}\Big]\Big)\nn\\ \label{subshd.13lk13}
	\eq
	where we defined $\alpha_k$ in the second inequality as follows:
	\begin{align}
	1-\alpha_k\define& \left(\mu_o \eta_k  -\mu_o^2 e_k^2  - \mu_o^2\beta_k^2 - 2\mu_o^2h e^2_k\right)\frac{q_k}{p_k}\nn\\
	\stackrel{\eqref{relation}}{=}& \mu_k\big(\eta_k - \mu_o e_k^2  - \mu_o\beta_k^2 - 2\mu_oh e^2_k\big)
	\end{align}
	Let $\alpha$ denote the largest $\alpha_k$ among all agents:
	\be
	\alpha\;\define \;\max_{1\leq k\leq N}\;\{\alpha_k\}.
	\ee
	Then, it holds that when $\alpha\in(0,1)$, which will be shown later in (\ref{see.dlkarqe}):
	\begin{align}
		\sum_{k=1}^N\left(1-\alpha_k\right) p_k\Ex\|\widetilde \w_{k,i-1}  \|^2
		\geq&\; (1-\alpha) \sum_{k=1}^N p_k\Ex\|\widetilde \w_{k,i-1}  \|^2\nn\\
		\geq&\; (1-\alpha)\Ex \| \bar \w_{i-1}\|^2,
	\end{align}
	where we used Jensen's inequality to deduce that
	\be
	\| \bar\w_{i-1}\|^2 =\left\|\sum_{k=1}^{N} p_k\widetilde\w_{k,i-1} \right\|^2 \leq  \sum_{k=1}^N p_k  \|\widetilde\w_{k,i-1}\|^2.
	\ee
	It follows from (\ref{subshd.13lk13}) that
	\begin{align}
	&2\mu_o \Big(\sum_{k=1}^{N} q_k \left(\Ex\,J_k(\w_{k,i-1}) - J_k(w^\star)\right)  \Big)\hspace{3.4cm}\nn\\
	&\;\;\;\leq \alpha\, \Ex \| \bar \w_{i-1}\|^2-\Ex \| \bar \w_{i}\|^2\nn\\
	&\hspace{7mm}{}+\mu_o^2 \sum_{k=1}^N\Big(q_k f_k^2+q_k\sigma^2_k + 2hq_k\Big[f_k^2+\|g'_k(w^\star)\|^2+\frac{1}{2}\Big]\Big)\label{9g.49}
	\end{align}
	This inequality recursion has a form similar to the one we encountered in the single agent case. Specifically, let us introduce the scalars:
	\begin{align}
		a(i) \define& \,\sum_{k=1}^{N} q_k \left(\Ex\,J_k(\w_{k,i-1}) - J_k(w^\star)\right) \\
		b(i) \define& \, \Ex \| \bar \w_{i}\|^2\\
		\tau^2 \define&\,\sum_{k=1}^N\Big(q_k f_k^2+q_k\sigma^2_k + 2hq_k\Big[f_k^2+\|g'_k(w^\star)\|^2+\frac{1}{2}\Big]\Big)
	\end{align}
	Then, recursion (\ref{9g.49}) can be rewritten more compactly in the form:
	\be
		2\mu_o a(i)\leq \alpha b(i-1) - b(i) + \mu_o^2 \tau^2 \label{e8g.23j}
	\ee
	This recursion has the same format as equation (69) in Part I \cite{ying16ssgd1}. Lastly, notice that 
	\begin{align}
		&\hspace{-6mm} \sum_{k=1}^N q_k \Big(\Ex J_k(\w_{k,i}^{\rm best}) - J_{k}(w^\star)\Big) \nn\\
		\stackrel{\eqref{w_k.best}}{=}&\, \sum_{k=1}^Nq_k \Big(\min_{1\leq i \leq L}\Ex J_k(\w_{k,i-1})- J_{k}(w^\star)\Big)\nn\\
		\leq&\,\min_{1\leq i \leq L}\sum_{k=1}^Nq_k \Big(\Ex J_k(\w_{k,i-1})- J_{k}(w^\star)\Big)\nn\\
		= &\, \min_{1\leq i \leq L} a(i) 
	\end{align}
	{\color{black} This result ensure that $\w^{\rm best}_{k,i}$ satisfies a condition similar to (76) in Part I \cite{ying16ssgd1}. 
	}
	The argument can now be continued similarly to arrive at the conclusions in the statement of the theorem.  Stability  is ensured by requiring $\alpha_k\in(0,1)$, i.e.,
	\begin{align}
	\alpha_k =&1- \mu_k\big(\eta_k - \mu_o e_k^2  - \mu_o\beta_k^2 - 2\mu_oh e^2_k\big)\in (0,1)
	\end{align}
	The condition $\alpha_k<1$ is met for 
	\begin{align}
	\mu_o < \frac{\eta_k}{\beta_k^2+(1+2h)e_k^2},\quad\forall k.
	\end{align}
	while the condition $\alpha_k>0$ requires
	\begin{align}
		\mu_k\big(\eta_k - \mu_o e_k^2  - \mu_o\beta_k^2 - 2\mu_oh e^2_k\big) < 1\label{see.fj8g}
	\end{align}
	But because $\eta_k - \mu_o e_k^2  - \mu_o\beta_k^2 - 2\mu_oh e^2_k\leq \eta_k$, we conclude $0<\mu_k< \frac{1}{\eta_k}$ is sufficient for condition (\ref{see.fj8g}). Combining these conditions with  (\ref{relation}), we establish
	\be
		\mu_k < \min\left\{\frac{1}{\eta_k}, \frac{\eta_k q_k}{p_k\beta^2_k +(1+2h)p_ke_k^2}\right\} \label{see.dlkarqe}
	\ee
	which ensures $\alpha_k\in (0,1)$.

	\section{\sc Proof of (\ref{prove.this.appendix})} \label{app.netClose}
	We establish the asymptotic result (\ref{prove.this.appendix}). Let
	\be
	\bar\swb_{i} = {\rm col}\{\bar\w_{i},\ldots,\bar\w_{i}\} = \one_N\otimes \bar\w_{i},
	\ee
	where the vector  $\bar\w_{i}$  is stacked $N$ times to match the dimension of $\widetilde{\swb}_i$. We start from the second relation in the error recursion (\ref{eq.errorReccursion2}):
	\be
	\check{\swb}_i = \mathcal{J}\tran_\epsilon\Big(\check{\swb}_{i-1} + \mu_o \check{\sgb}(\swb_{i-1})+\mu_o \check{\ssb}_i(\swb_{i-1})\Big)\label{eq.eq154},
	\ee
	{\color{black} and first explain how to recover $\widetilde\swb_{i} - \bar\swb_{i}$ from $\check{\swb}_i$. From (\ref{def.w2})
	\begin{align}
	 \widetilde\swb_{i} =&\;\cV_{\epsilon}^{-\mathsf{T}} \ba{c}\bar\w_{i}\\ \check\swb_{i}\ea\nn\\
		\stackrel{(\ref{eig.decomp})}{=}&\;\ba{c|c}\one\otimes I_M&\cV_L\ea \ba{c}\bar\w_{i}\\ \check\swb_{i}\ea\nn\\
		=&\; \bar\swb_{i} + \cV_L\check\swb_{i} \label{fiu29jf.e3}
	\end{align}
	}
	Next, returning to the error recursion (\ref{eq.eq154}), and computing the expected squared norm, we obtain:
	\bq
	\Ex[\|\check{\swb}_i\|^2\,|\,{\bm{\cal F}}_{i-1}] \hspace{-2mm}&=&\hspace{-2mm} \left\| \mathcal{J}\tran_\epsilon\Big(\check{\swb}_{i-1} + \mu_o  \check{\sgb}(\swb_{i-1})\Big)\right\|^2\nn\\
	&&{}+\mu_o^2\Ex[ \| \mathcal{J}\tran_\epsilon \check{\ssb}_i(\swb_{i-1})\|^2 \,|\,{\bm{\cal F}}_{i-1}] \nn\\
	&\leq&\hspace{-2mm} \rho(\mathcal{J}_\epsilon\mathcal{J}\tran_\epsilon)\| \check{\swb}_{i-1} + \mu_o  \check{\sgb}(\swb_{i-1})\|^2\nn\\
	&& +\mu_o^2\rho(\mathcal{J}_\epsilon\mathcal{J}\tran_\epsilon)\Ex[ \| \check{\ssb}_i(\swb_{i-1})\|^2 \,|\,{\bm{\cal F}}_{i-1}], \nn\\\label{eq.diffRecursion}
	\eq
	where, from \cite[Ch.~9]{sayed2014adaptation}, we know that
	\be
	\rho(\mathcal{J}_\epsilon\mathcal{J}\tran_\epsilon) \leq (\rho(J_\epsilon)+\epsilon)^2<1.	
	\ee
	Let us examine the terms in (\ref{eq.diffRecursion}). To begin with, note that
	\bq
	&&\hspace{-9mm}\rho(\mathcal{J}_\epsilon\mathcal{J}\tran_\epsilon)\| \check{\swb}_{i-1} + \mu_o  \check{\sgb}(\swb_{i-1})\|^2   \label{eq.eq162}\nn\\
	&\leq&\hspace{-2mm} (\rho(J_\epsilon)+\epsilon)^2 \left\|t\frac{1}{t} \check{\swb}_{i-1} +\frac{1-t}{1-t} \mu_o  \check{\sgb}(\swb_{i-1})\right\|^2\nn\\
	&\stackrel{(a)}{\leq}& \hspace{-2mm}\frac{(\rho(J_\epsilon)+\epsilon)^2}{t} \| \check{\swb}_{i-1} \|^2+ \mu_o^2\frac{(\rho(J_\epsilon)+\epsilon)^2}{1-t}  \|\check{\sgb}(\swb_{i-1})\|^2\nn\\
	&\stackrel{(b)}{\leq}&\hspace{-2mm} (\rho(J_\epsilon)+\epsilon) \| \check{\swb}_{i-1} \|^2+ \mu_o^2\frac{(\rho(J_\epsilon)+\epsilon)^2}{1-\rho(J_\epsilon)-\epsilon}  \|\check{\sgb}(\swb_{i-1})\|^2,\nn\\
	\eq
	where step (a) is because of Jensen's inequality and in step (b) we select $t=\rho(J_\epsilon)+\epsilon<1$.
	Next, we bound the square of the sub-gradient term:
	\begin{align}
	&\hspace{-3mm}\|\check{\sgb}(\swb_{i-1})\|^2 \nn\\
	{=}&\, \|\mathcal{V}_R\tran\cU \sgb(\swb_{i-1})\|^2\nn\\
	\leq&\,  \|V_R\|^2\|U\|^2\Big(\sum_{k=1}^N\|g_k(\w_{k,i-1})\|^2\Big) \nn\\
	\stackrel{(a)}{\leq}&\,2  \|V_R\|^2\|U\|^2\Big(\sum_{k=1}^N\|g_k(\w_{k,i-1}) - g'_k(w^\star)\|^2 +\|g'_k(w^\star)\|^2\Big) \nn\\
	\leq&\,2 \|V_R\|^2\|U\|^2\Big(\sum_{k=1}^Ne^2_k\|\widetilde\w_{k,i-1}\|^2+f^2_k+\|g'_k(w^\star)\|^2\Big)\nn\\
	\stackrel{(b)}{\leq}&\,2 \|V_R\|^2\|U\|^2\Big(e^2_{\max}\|\widetilde\swb_{i-1}\|^2+\sum_{k=1}^N(f^2_k+\|g'_k(w^\star)\|^2)\Big),
	\end{align}
	where in step (a) we subtract and add $g'_k(w^\star)$ inside of the norm and the factor 2 comes from Jensen's inequality, and in step (b) we let $e^2_{\max}=\max_k e_k^2 $. We can then bound (\ref{eq.eq162}) by
	\bq
	&&\hspace{-10mm}\rho(\mathcal{J}_\epsilon\mathcal{J}\tran_\epsilon)\| \check{\swb}_{i-1} + \mu  \check{\sgb}(\swb_{i-1})\|^2 \nn\\ 
	&\hspace{-4mm}\leq&\hspace{-2mm} (\rho(J_\epsilon)+\epsilon) \| \check{\swb}_{i-1} \|^2\nn\\
	&&\hspace{-1mm} {}+ 2\mu^2_o\frac{(\rho(J_\epsilon)+\epsilon)^2}{1-\rho(J_\epsilon)-\epsilon}  \|V_R\|^2\|U\|^2e^2_{\max}\|\widetilde\swb_{i-1}\|^2
	\nn\\
	&&\hspace{-1mm} {}+2\mu^2_o\frac{(\rho(J_\epsilon)+\epsilon)^2}{1-\rho(J_\epsilon)-\epsilon}  \|V_R\|^2\|U\|^2\sum_{k=1}^N (f^2_k+\|g'_k(w^\star)\|^2).\nn\\
	\eq
	Finally, we consider the last term involving the gradient noise in (\ref{eq.diffRecursion}):
	\bq
	&&\hspace{-12mm}
	\Ex[ \| \check{\ssb}_i(\swb_{i-1})\|^2 \,|\,{\bm{\cal F}}_{i-1}]\nn\\
	&\leq&\hspace{-2mm} \|V_R\|^2\|U\|^2\left(\sum_{k=1}^N \beta_k^2\|\widetilde\w_{k,i-1}\|^2+\sigma_k^2 \right) \nn\\
	&\leq&\hspace{-2mm} \|V_R\|^2\|U\|^2 \beta_{\max}\|\widetilde{\swb}_{i-1}\|^2 + \|V_R\|^2\|U\|^2\sum_{k=1}^N\sigma_k^2.
	\eq
	Now introduce the constants:
	\bq
	a&\hspace{-0.3cm}\define\hspace{-0.3cm}&\frac{2(\rho(J_\epsilon)	+\epsilon)^2}{1-\rho(J_\epsilon)-\epsilon}  \|V_R\|^2\|U\|^2e^2_{\max}\nn\\
	&&\;\;{}+\,\rho(J_\epsilon  J\tran_\epsilon)\|V_R\|^2\|U\|^2 \beta_{\max},\\
	b&\hspace{-0.3cm}\define\hspace{-0.3cm}& \frac{2(\rho(J_\epsilon)+\epsilon)^2}{1-\rho(J_\epsilon)-\epsilon}\|V_R\|^2\|U\|^2\sum_{k=1}^N (f^2_k+\|g'_k(w^\star)\|^2)\nn\\
	&&\;\; {}+\rho(J_\epsilon  J\tran_\epsilon)\|V_R\|^2\|U\|^2\sum_{k=1}^N\sigma_k^2.
	\eq
	Although the matrix $U$ is dependent on the $\mu_k$, entries of $U$ are ratios relative to $\mu_o$.
	Then, substituting the previous results into  (\ref{eq.diffRecursion}), we arrive at
	\be
	\Ex\| \check{\swb}_{i}\|^2 \leq (\rho(J_\epsilon)+\epsilon)\Ex \| \check{\swb}_{i-1}\|^2 + \mu_o^2a\Ex\|\widetilde{\swb}_{i-1}\|^2 + \mu_o^2b \label{eq.160},
	\ee
	In Appendix \ref{app.c} we show that $\Ex\|\widetilde{\swb}_{i-1}\|^2$, for any iteration $i$, is bounded by a constant value for sufficient small step-sizes. In this case, we can conclude that
	\be
	\Ex\| \check{\swb}_{i}\|^2 \leq (\rho(J_\epsilon)+\epsilon)\Ex \| \check{\swb}_{i-1}\|^2 + \mu_o^2b',
	\ee
	for some constant $b'$, so that at steady state:
	\be
	\limsup_{i\to\infty}\Ex\| \check{\swb}_{i}\|^2 \leq \frac{\mu_o^2 b'}{1-\rho(\mathcal{J}_\epsilon)-\epsilon}=O(\mu_o^2)\label{eq.eq177}.
	\ee
	Using relation (\ref{fiu29jf.e3}), it then follows asymptotically that
	for $i\gg~1$:
	\bq
	\Ex\|\widetilde{\swb}_i - \bar{\swb}_i\|^2&\leq&  \|V_L\|^2\cdot \Ex\| \check{\swb}_{i}\|^2 =O(\mu_o^2),
	\eq
	and, consequently,
	\bq
	\Ex\|\widetilde{\w}_{k,i} - \bar{\w}_{i}\|^2&\leq& \Ex\|\widetilde{\swb}_i - \bar{\swb}_i\|^2=O(\mu_o^2).
	\eq
	
	\section{\sc Proof that $\Ex\|\widetilde{\swb}_{i}\|^2$ is uniformly bounded} \label{app.c}
	We follow mathematical induction to establish that $\Ex\|\widetilde{\swb}_i\|^2$ is uniformly bounded by a constant value, for all $i$. Assume, at the initial time instant we have $\Ex\|\widetilde{\w}_{k,0}\|^2<\mathsf{c}$ for all $k$ and for some constant value $\mathsf{c}$.
	Then, assuming this bound holds at iteration $i-1$, namely,
	\be
	\Ex\|\widetilde{\w}_{k,i-1}\|^2 \leq \mathsf{c},\quad\forall k\label{eq.166},
	\ee
	we would like to show that it also holds at iteration $i$. Recall from (\ref{atc.eq2}) that the diffusion strategy consists of two steps: an adaptation step followed by a combination step. The adaptation step has a similar structure to the single-agent case. Hence, the same derivation that was used to  establish for single agent case in Part I\cite[Eq. 64]{ying16ssgd1} would show that for agent $k$:
	\bq
	&&\hspace{-12mm}2\mu_k \left(\Ex\,J_k(\w_{k,i-1}) - J_k(w^\star_k)\right) \nn\\
	&\leq&\hspace{-2mm} \alpha_k \Ex \| \widetilde\w_{k,i-1}\|^2 
		{}-\Ex \| \widetilde\bpsi_{k,i}\|^2  + \mu_k^2 (f_k^2+\sigma^2_k),
	\eq
	where
	\begin{align}
	\alpha_k =\,& 1-\mu_k\eta_k+\mu_k^2(e_k^2+\beta_k^2) = 1-O(\mu_k),\\
	w^\star_k \define& \argmin_w  J_k(w).
	\end{align}
	Now, since $\Ex\,J_k(\w_{k,i-1})\geq J_k(w^\star_k)$, we conclude that
	\begin{align}
	\Ex \| \widetilde\bpsi_{k,i}\|^2&\;\,\leq\;\; \alpha_k \Ex \| \widetilde\w_{k,i-1}\|^2+ \mu_k^2 (f_k^2+\sigma^2_k)\nn\\
	&\stackrel{(\ref{eq.166})}{\leq}\alpha_k \mathsf{c}+ \mu_k^2 (f_k^2+\sigma^2_k),
	\end{align}
	where the step-size $\mu_k$ can be chosen small enough to ensure $\alpha_k\in (0,1)$. Now, it is also clear that
	there exist sufficiently small values for $\mu_k$ to ensure that, for all agents
	$k$:
	\be
	\alpha_k \mathsf{c} +  \mu_k^2 (f_k^2+\sigma^2_k) \leq \mathsf{c},
	\ee
	which then guarantees that
	\be
	\Ex \| \widetilde\bpsi_{k,i}\|^2\leq \mathsf{c}.
	\ee
	It then follows from the combination step (\ref{atc.eq2}) that
	\bq
	\Ex \| \widetilde\w_{k,i}\|^2 &=& \Ex \left\|\sum_{\ell\in\mathcal{N}_k} a_{\ell k}\widetilde\bpsi_{\ell,i}\right\|^2\nn\\
	&\leq&\sum_{\ell\in\mathcal{N}_k} a_{\ell k}\Ex\left\|\widetilde\bpsi_{\ell,i}\right\|^2\nn\\
	&\leq& \sum_{\ell\in\mathcal{N}_k} a_{k \ell }\mathsf{c}\nn\\
	&=& \mathsf{c},\;\;\forall k.
	\eq
	Therefore, starting from (\ref{eq.166}), we conclude that
	$\Ex\|\widetilde{\w}_{k,i}\|^2< \mathsf{c}$ as well, as desired. Finally, since $\Ex\|\widetilde{\swb}_i\|^2 = \sum_{k=1}^N \Ex \| \widetilde\w_{k,i}\|^2$, we conclude that $\Ex\|\widetilde{\swb}_i\|^2$ is also uniformly bounded over time.
	
	\section{Proof of corollary \ref{corollary.1}} \label{app.corollary1}
	Iterating (\ref{e8g.23j}) over $1\leq i\leq L$, for some interval length $L$, gives:
	\be
	\sum_{i=1}^L\alpha_m^{L-i} (2\mu_oa(i) - \mu_o^2 \tau^2)\leq \alpha_m^Lb(0)	
	\ee
	Then, dividing both side by the same sum:
	\be
	\sum_{i=1}^L\frac{\alpha_m^{L-i}}{S_{L-1}} (2\mu_oa(i) - \mu_o^2 \tau^2)\leq \frac{\alpha_m^L}{S_{L-1}}b(0)	\label{2g89gf}
	\ee
	Now, because of the convexity of each $J_k(\cdot)$, we have
	\be
	J_k(\bar{\w}_{k,{L-1}})\leq\sum_{j=0}^{L-1} r_{L-1}(j)J_k(\w_{k,j})
	\ee	
	Thus, we can establish:
	\begin{align}
	&\hspace{-4mm}\sum_{i=1}^L\frac{\alpha_m^{L-i}}{S_{L-1}} a(i)\nn\\
	=&\, \sum_{i=1}^L r_{L-1}(i-1) \sum_{k=1}^Nq_k\Big(\Ex J_k(\w_{k,i-1})-J_k(w^\star) \Big)\nn\\
	\geq&\,\sum_{k=1}^N q_k \Big(\Ex J_k(\bar{\w}_{k,L-1}) - J_k(w^\star)\Big)
	\end{align}
	Substituting into (\ref{2g89gf}), we establish:
	\be
	2\mu_o\sum_{k=1}^N q_k \Big(\Ex J_k(\bar{\w}_{k,L-1}) - J_k(w^\star)\Big) \leq \frac{\alpha_m^L}{S_{L-1}}b(0)+\mu_o^2\tau^2
	\ee
	Letting $L\to\infty$, we establish \eqref{llkad.lk1l3kmjkda}.
	
\bibliographystyle{IEEEbib}
\bibliography{ref_sbg}
\end{document}